\documentclass[11pt]{amsart}
\usepackage{amssymb,amsmath,mathrsfs,graphics,graphicx}
\usepackage[T1]{fontenc}
\usepackage{float}

\theoremstyle{plain}
\newtheorem{thm}{Th\'eor\`eme}
\newtheorem{cor}[thm]{Corollaire}

\theoremstyle{definition}

\def\vs{\vspace{0.3cm}}
\def\hs{\hspace{0.1cm}}

\def\C{\mathbb{C}}
\def\P{\mathbb{P}}
\def\Z{\mathbb{Z}}

\def\Bir{\mathsf{Bir}}

\def\PGL{\mathsf{PGL}}
\def\SL{\mathsf{SL}}

\newcommand{\transposee}[1]{{\vphantom{#1}}^{\mathit t}{#1}}

\def\og{\leavevmode\raise.3ex\hbox{$\scriptscriptstyle\langle\!\langle$~}}
\def\fg{\leavevmode\raise.3ex\hbox{~$\!\scriptscriptstyle\,\rangle\!\rangle$}}

\addtocounter{section}{0}             % Start with section 1
\numberwithin{equation}{section}       % Number formulas within sections
\begin{document}

\title{Le groupe de \textsc{Cremona} est hopfien}
\author{Julie D\'eserti}
\email{julie.deserti@univ-rennes1.fr}
\maketitle

\begin{abstract}
On d\'ecrit les endomorphismes du groupe de \textsc{Cremona}
et on en d\'eduit son caract\`ere hopfien.
\end{abstract}

\vspace{3mm}

\noindent Une transformation rationnelle de $\P^2(\C)$ dans 
lui-m\^eme s'\'ecrit $$(x:y:z)\mapsto
(P_0(x,y,z):P_1(x,y,z):P_2(x,y,z))$$ o\`u les $P_i$ d\'esignent
des polyn\^omes homog\`enes de m\^eme degr\'e. Lorsqu'elle 
est inversible, on dit qu'elle est birationnelle~; par exemple 
l'involution de \textsc{Cremona} $\sigma=(yz:xz:xy)$ est 
birationnelle. Le groupe des transformations birationnelles, not\'e 
$\Bir(\P^2(\C)),$ est aussi appel\'e groupe de \textsc{Cremona}.

\begin{thm}[N\oe ther, \cite{AC,Ca}]\label{nono}
{\sl Le groupe de \textsc{Cremona} est engendr\'e par 
$\PGL_3(\C)$ et l'invo\-lution $\sigma=(yz:xz:xy).$}
\end{thm}

\noindent Un automorphisme
$\tau$ du corps $\C$ induit un isomorphisme $\tau(.)$ de 
$\Bir(\P^2(\C))~:$ \`a un \'el\'ement $f$ de $\Bir(\P^2(\C))$ 
nous associons l'\'el\'ement
$\tau(f)$ obtenu en faisant agir $\tau$ sur les coefficients de
$f$ exprim\'e en coordonn\'ees homog\`enes. {\sl Tout 
automorphisme du groupe de \textsc{Cremona} s'obtient \`a
partir de l'action d'un automorphisme de corps et d'une 
conjugaison int\'erieure (\cite{De2}).}

\noindent Ici nous nous int\'eressons aux endomorphismes du 
groupe de \textsc{Cremona}~:

\begin{thm}\label{end}
{\sl Soit $\varphi$ un endomorphisme non trivial de $\Bir(
\P^2(\C)).$ Il existe un plongement de corps $\lambda$ de $\C$ dans 
lui-m\^eme et une transformation birationnelle $\psi$ tels
que pour tout $f$ dans $\Bir(\P^2(\C))$ on ait $$\varphi(f)=
\lambda(\psi f\psi^{-1}).$$ En particulier $\varphi$ est 
injectif.}
\end{thm}

\noindent Une cons\'equence directe est la suivante~:

\begin{cor}
{\sl Le groupe de \textsc{Cremona} est hopfien, {\it i.e.} tout
endomorphisme surjectif de $\Bir(\P^2(\C))$ est un 
automorphisme.}
\end{cor}

%\begin{proof}[{\sl D\'emonstration}]
%Soit $\varphi$ un endomorphisme surjectif de $\Bir(\P^2(\C))$
%dans lui-m\^eme. D'apr\`es le Th\'eor\`eme \ref{end} il 
%existe une immersion $\lambda$ de $\C$ dans lui-m\^eme telle que
%$\varphi(\Bir(\P^2(\C)))=\Bir(\P^2(\lambda(\C)))~;$ mais $\varphi$
%\'etant surjectif, on a $\Bir(\P^2(\C))=\Bir(\P^2(\lambda(\C))).$ 
%Par suite $\C=\lambda(\C)$, {\it i.e.} $\lambda$ est un 
%automorphisme de corps. Le Th\'eor\`eme \ref{end} assure que 
%$\varphi$ est un automorphisme du groupe de \textsc{Cremona}.
%\end{proof}

\noindent La preuve du Th\'eor\`eme \ref{end} repose
en partie sur le r\'esultat suivant que nous 
appliquons \`a $\Gamma=~\SL_3(\Z)~:$

\begin{thm}[\cite{De2}]\label{rigib}
{\sl Soient $\Gamma$ un sous-groupe d'indice fini de
$\SL_3(\Z)$ et $\rho$ un morphisme injectif de $\Gamma$ dans
$\Bir(\P^2(\C)).$ Alors $\rho$ co\"incide, \`a conjugaison
birationnelle pr\`es, avec le plongement
canonique ou la contra\-gr\'ediente, {\it i.e.} l'involution
$u\mapsto \transposee u^{-1}$.}
\end{thm}

\vs

\noindent On travaille dans une carte affine $(x,y)$ de 
$\P^2(\C).$ Introduisons le groupe des translations~:  
$$\mathsf{T}=\{(x+\alpha,y+\beta)\hs|\hs\alpha,\hs\beta\in
\C\}.$$ %Le groupe $\tilde{\mathsf{T}}$ d\'esignera 
%$$\left\{\left(\frac{x}{\alpha x+\beta y+1},\frac{y}{\alpha 
%x+\beta y+1}\right)\hs|\hs\alpha,\hs\beta\in\C\right\}~;$$ 
%il s'agit de l'image de $\mathsf{T}$ par la contragr\'ediente.

\begin{proof}[{\sl D\'emonstration du Th\'eor\`eme \ref{end}}]
\noindent Puisque $\PGL_3(\C)$ est simple,
$\varphi_{|\PGL_3(\C)}$ est ou bien triviale, ou bien injective.\vs

\noindent{\sl 1.} Supposons $\varphi_{|\PGL_3(\C)}$ triviale. Posons 
$h:=~(x,x-y,x-z)$~; comme l'a remarqu\'e \textsc{Gizatullin} 
(\cite{Gi}), on a $(h\sigma)^3=\mathsf{id}.$
Ainsi $\varphi((h\sigma)^3)=\varphi(\sigma)=
\mathsf{id},$ {\it i.e.} $\varphi$ est trivial d'apr\`es le 
Th\'eor\`eme \ref{nono}.\vs

\noindent{\sl 2.} Si $\varphi_{|\PGL_3(\C)}$ est injective, alors
$\varphi_{|\SL_3(\Z)}$ est, \`a conjugaison birationnelle
pr\`es, le plongement canonique ou la contragr\'ediente.

\noindent{\sl 2.a.} Supposons que $\varphi_{|\SL_3(\Z)}=
\mathsf{id}.$ Notons $\mathsf{H}$ le
groupe des matrices $3\times 3$ triangulaires sup\'erieures 
unipotentes.
%$$\mathsf{H}:=\left\{\left(
%\begin{array}{ccc}
%1 & a & b\\
%0 & 1 & c\\
%0 & 0 & 1
%\end{array}
%\right)\hs|\hs a,\hs b,\hs c\in\C\right\}.$$ 
Posons~: $$f_\beta(x,y):=\varphi(x+\beta,y),\hs \hs g_\alpha(x,
y):=\varphi(x+\alpha y,y)\hspace{3mm} \text{et}\hspace{3mm} h_\gamma(x,y):=\varphi(x,y+\gamma).$$ Les transformations
birationnelles $f_\beta$ et $h_\gamma$ commutent \`a $(x+1,y)$ 
et $(x,y+1)$ donc $$f_\beta=(x+\lambda(\beta),y+\zeta(\beta))
\hspace{6mm}\text{et}\hspace{6mm}h_\gamma=(x+\eta(\gamma),y+
\mu(\gamma))$$ o\`u $\eta,\hs\zeta,\hs\mu$ et $\lambda$ sont 
des morphismes additifs de $\C$ dans $\C~;$ puisque $g_\alpha$ commute 
\`a $(x+y,y)$ et $(x+1,y)$ il est de la forme $(x+A_\alpha(y),
y)$. La relation $$(x+\alpha y,y)(x,y+\gamma)(x+\alpha y,
y)^{-1}(x,y+\gamma)^{-1}=(x+\alpha\gamma,y)$$ implique
que, pour tous nombres complexes $\alpha$ et $\gamma$, nous avons
$g_\alpha h_\gamma=f_{\alpha\gamma}h_\gamma g_\alpha$. Nous 
en d\'eduisons que~: $$f_\beta=(x+\lambda(\beta),y),
\hspace{6mm}g_\alpha=(x+\Theta(\alpha)y+\varsigma(\alpha), y)
\hspace{6mm}\text{et}\hspace{6mm}\Theta(\alpha)\mu(\gamma)=
\lambda(\alpha\gamma).$$ En utilisant l'\'egalit\'e $$(x+
\alpha)(x,\beta x+y)(x-\alpha,y)(x,y-\beta x)=(x,y-\alpha
\beta)$$ on \'etablit que $h_\gamma=(x,y+\mu(\gamma)).$ Autrement
dit $$\varphi(x+\alpha,y+\beta)=(x+\lambda(\alpha),y+\mu(\beta))
\hspace{6mm}\forall\hs\alpha,\hs\beta\in\C.$$ Ainsi 
$\varphi(\mathsf{T})\subset\mathsf{T}$ et $\varphi(\mathsf{H})
\subset\mathsf{H}~;$ puisque $\PGL_3(\C)$ est engendr\'e par 
$\mathsf{H}$ et $\SL_3(\Z),$ l'image de $\PGL_3(\C)$ par
$\varphi$ est contenue dans $\PGL_3(\C).$ 
Le Th\'eor\`eme de classification de \textsc{Borel} et 
\textsc{Tits} (\cite{BoTi}, Th\'eor\`eme A, p. 500)
assure qu'\`a conjugaison int\'erieure pr\`es  
l'action de $\varphi$ sur $\PGL_3(\C)$ provient d'un
plongement de corps de $\C$ dans lui-m\^eme.

\noindent{\sl 2.b.} Supposons que la restriction de $\varphi$ 
\`a $\SL_3(\Z)$ co\"{\i}ncide avec la contragr\'ediente. 
En \'etudiant 
les images de $\mathsf{T}$ et $\mathsf{H}$ par $\varphi$, on 
montre que $\varphi(\PGL_3(\C))\subset\PGL_3(\C).$ 
Toujours d'apr\`es \cite{BoTi} (Th\'eor\`eme A, p. 500) \`a 
conjugaison int\'erieure 
pr\`es, l'action de $\varphi$ sur $\PGL_3(\C)$ provient ici 
d'un plongement de corps de $\C$ dans lui-m\^eme compos\'e avec la 
contragr\'ediente.

\vs

\noindent{\sl 3.} Supposons donc que l'action de $\varphi$
sur $\PGL_3(\C)$ co\"{\i}ncide avec celle d'un plongement de
corps
$\lambda$ de $\C$ dans lui-m\^eme ou avec la compos\'ee d'une 
telle action et de la contragr\'ediente. 

\noindent Posons $(\tau_1,\tau_2)=\varphi(x,1/y).$
\`A partir de $$(x,1/y)(\alpha x,\beta y)(x,1/y)=
(\alpha x,y/\beta)$$ on obtient $$\tau_1(\lambda(
\alpha^{-1})x,\lambda(\beta^{-1})y)=\lambda(\alpha^{-1})
\tau_1(x,y)\hspace{6mm}\text{et}\hspace{6mm}\tau_2(\lambda(
\alpha^{-1})x,\lambda(\beta^{-1})y)=\lambda(\beta)
\tau_2(x,y)$$ ou $$\tau_1(\lambda(\alpha)x,
\lambda(\beta)y)=\lambda(\alpha)\tau_1(x,y)\hspace{6mm}\text{et}
\hspace{6mm}\tau_2(\lambda(\alpha)x,\lambda(\beta)y)=\frac{
\tau_2(x,y)}{\lambda(\beta)}$$ suivant que la contragr\'ediente
intervient ou non. Par suite $\varphi(x,1/y)=(\pm x,\pm 1/y).$ 

\noindent L'\'egalit\'e $((y,x)(x,1/y))^2=\sigma$ assure que $\varphi(\sigma)=\pm\sigma.$ Notons $h:=\left(\frac{
x}{x-1},\frac{x-y}{x-1}\right)~;$ la transformation $(h
\sigma)^3$ est triviale (\emph{voir} \cite{Gi}) donc $(\varphi(h)\varphi(\sigma))^3$ doit aussi 
l'\^etre. Puisque $h$ appartient \`a $\SL_3(\Z),$ on a 
$\varphi(h)=h$ ou $\varphi(h)=(-x-y-1,y)$ suivant que 
$\varphi_{|\SL_3(\Z)}$ est l'identit\'e
ou la contragr\'ediente. Si $\varphi(h)=h$, alors 
$\varphi(\sigma)=\sigma$ et on conclut avec le Th\'eor\`eme
\ref{nono}.
Lorsque $\varphi(h)=(-x-y-1,y),$ la seconde composante de
$(\varphi(h)\varphi(\sigma))^3$ vaut $\pm 1/y$ ce qui est
exclu.
\end{proof}

\noindent{\bf Remerciements.}  Je remercie \'E. \textsc{Ghys} 
d'avoir pos\'e le probl\`eme ainsi que pour les remarques et 
suggestions qu'il m'a faites. Merci \`a D. \textsc{Cerveau} 
pour nos discussions anim\'ees et fructueuses.

\end{document}